\documentclass[12pt]{article}
\usepackage{cite}
\usepackage{amsmath,amssymb,amsfonts}
\usepackage{algorithmic}
\usepackage{graphicx}
\usepackage{textcomp}

\usepackage{epsfig}
\usepackage{latexsym}
\usepackage{pdfsync}

\def\Real{{I\!\!R}}

\newcommand{\eq}{\begin{equation}\begin{array}{rllllllllllllllllllllllllllllllll}}
\newcommand{\ee}{\end{array}\end{equation}}
\newcommand{\bmt}{\left[ \begin{array}{ccccccccc}}
\newcommand{\emt}{\end{array}\right]}
\newcommand{\bea}{\begin{eqnarray}}
\newcommand{\eea}{\end{eqnarray}}
\newcommand{\bean}{\begin{eqnarray*}}
\newcommand{\eean}{\end{eqnarray*}}
\newcommand{\bc}{\begin{center}}
\newcommand{\ec}{\end{center}}

\begin{document}
\title{Linear Quadratic Gaussian Synthesis for a Heated/Cooled  Rod Using Point Actuation and Point Sensing }
\thanks{This work was supported by AFOSR under FA9550-20-1-0318 }

\author{Arthur J. Krener\\
Naval Postgraduate School \\Monterey CA 93943-5216\\
\tt{ajkrener@nps.edu}}
\date{}

\maketitle

\begin{abstract}
We consider a rod  that is heated/cooled and sensed  at multiple point locations.   To stabilize it to a constant temperature
we set up  a Linear Quadratic Regulator that we explicitly solve by the method of completing the square to find the optimal linear
state feedback for the point actuators.  But we don't assume that the whole state is mesureable so we construct an infinite dimensional  Kalman filter to estimate
the whole state from a finite number of noisy point measurements.   These two components yield a Linear Quadratic Gaussian (LQG)  Synthesis for the heat equation
under point actuation and point sensing.
\end{abstract}

\section{Introduction}
  Linear Quadratic Gaussian Synthesis  is a (perhaps the) standard approach to constructing a compensator  for a finite 
  or an infinite dimensional linear system.   It consists  of two parts.  One first solves a Linear Quadratic Regulator problem for
  a linear state  feeback law for a finite number of actuators. If the state were fully measureable this linear state feedback would asymptotically stabilize the system.  But the state is
  usually not fully measureable, instead a finite number of linear functionals of state are available and these measurements are corrupted by noise.
  A Kalman filter is used to process the measurents to get an estimate of the whole state.  The certainty equivalence principle is envoked and  the linear feedback is applied to estimate of the state.  The result is a compensator of the same dimension as the original system.  It can be shown that
  spectrum of the combined system and its compensator is the union of the spectrum of system under linear full state feedback with the spectrum of the error
  dynamics of the Kalman filter.  If both these spectra lie in the open left half plane then the compensator will stabilize the original system.
  
  This LQG approach works over both finite and infinite time horizons.  In this paper we will only treat infinite  horizons.  To find the gain of the optimal
  linear  state feedback one has to solve an algebraic Riccati equation and to find the gain of the Kalman filter one has to solve a dual algebraic
  Riccati equation.   Good software exists to solve these equations in low to medium dimensions but they can be difficult to solve in high or infinte
  dimensions.   In infinite dimensions other difficulties can arise.  The actuation could be at points in the spatial domain, e.g., boundary control.  
  The measuements could also be at points.   Such systems are not "state linear systems" in the definition of Curtain an Zwart \cite{CZ95}, \cite{CZ20}.  State linear systems must
  have bounded input and output linear functionals. 
  The usual approach to dealing with a system with point actuation and/or point sensing is to approximate it by a state linear system.   Boundary control actuation is replaced
 by intense actuation over a short interval adjacent to the boundary,  in effect the control input multiplies a shaping function that approximates a delta function at the boundary.  Point sensing is replaced
by integration of the state against a shaping function that approximates a delta function at the measurment point.   For more details about this we refer the reader to Chapter 6 of \cite{CZ95} and Chapter 9 of \cite{CZ20} and their extensive references.   For more on boundary control of systems described by PDEs see the treatises of Lions \cite{JL71}, Lasiecka-Triggiani \cite{LT00} and Krstic-Smyshlyaev \cite{KS08}.  Hulsing \cite{Hu99} and Burns-Hulsing have addressed the computational issues associated with boundary control.

We take a different approach, in effect, we model boundary and other point actuators by delta functions and we model point sensing also by delta functions.  We use a method that we call completing the square
to overcome the mathematical  technicalities associated with these  delta functions.  To keep the discussion concrete we limit our consideration to a rod heated/cooled at boundary and other points.  We make noisy meaurements of its 
temperature at some other points.  Because we focus on systems modeled by the heat equation we are able to give explicit solutions to the LQR and Kalman filtering equations using the simple technique of completing the square.
In particular the infinte dimensional analogs of the algebraic Riccati equations are elliptic PDEs that we call Riccati PDEs.   These Riccati PDEs can be explicitly  solved in terms of the eigenfunctions of the Laplacian.  We restrict our attention 
to Neumann boundary conditions but our methods readily extend to other self adjoint boundary conditions for the Laplacian. 
We started using the completing the square technique on a distributed control problem \cite{Kr20a} but we quickly realized that it works well for boundary control problems.
 In \cite{Kr20b} we treated the  LQR control of a rod heated/cooled at one end and insulated at the other.  The boundary conditions were Neumann at the insulated end and Robin at the controlled end.  We have also used the completing the square technique for the LQR boundary control of the wave equation \cite{Kr21a} and the beam equation \cite{Kr21b}.

 The rest of the paper is as follows.  In the next section we treat the LQR control of the rod heated/cooled at the boundary and other points.  Section 3 contains an example of heating/cooling at both ends of the rod and Section 4 contains an example of heating/cooling at both ends of the rod and at its midpoint.   In Section 5 we derive the Kalman filter for the rod with noisy measurements at several points by
 converting the filtering problem into a family of LQR problems.  Section 6 contains an example of the LQG synthesis of a Kalman filter with an LQR state feedback law.  The conclusion is found in Section 7.

\section{LQR for a Rod Heated/Cooled at Multiple Points}
Consider a rod of length one that is heated/cooled at multiple locations.  We let $x\in [0,1]$ denote position
on the rod  and let $z(x,t)$ denote the temperature of the rod at position $x$ at time $t$.  We assume the 
rod can be heated or cooled at $0=\xi_1 <\xi_2<\ldots< \xi_m=1$.  We let $u_k(t)$ be the heating/cooling flux applied to the
rod at $\xi_k$ for $k=1,\ldots,m$ and $u(t)=[u_1(t), \ldots u_m(t)]'$.  We model the rod by these equations
\bean
0&=&\frac{\partial z}{\partial t} (x,t) -\frac{\partial^2 z}{\partial x^2}(x,t) \\
0&=&z(\xi_k^-,t) -z(\xi_k^+,t) \mbox{ for } k=2,\ldots, m-1\\
\beta_k u_k(t)&=&\frac{\partial z}{\partial x} (\xi_k^-,t) -\frac{\partial z}{\partial x} (\xi_k^+,t) \mbox{ for } k=1,\ldots, m
\\
z(x,t) &=&  z^0(x)
\eean
where 
\bean
z(\xi_k^+,t)=\lim_{x\to \xi_k^+}z(x,t),&&z(\xi_k^-,t)=\lim_{x\to \xi_k^-}z(x,t)\\
\frac{\partial z}{\partial x} (\xi_k^+,t) =\lim_{x\to \xi_k^+}\frac{\partial z}{\partial x} (x,t),&&\frac{\partial z}{\partial x} (\xi_k^-,t) =\lim_{x\to \xi_k^-}\frac{\partial z}{\partial x} (x,t)
\eean
with the understanding that $\frac{\partial z}{\partial x} (x,t) =0$ outside of $[0,1]$ and $\beta_k\ge 0$.   If the rod is not heated/cooled at its endpoints we set $\beta_1=0$ and/or $\beta_m=0$.

The open loop system, $ u_k(t)=0$ for $k=1,\ldots,m$,  reduces to the standard heat equation with Neumann boundary conditions,
\bean
\frac{\partial z}{\partial t} (x,t) &=&\frac{\partial^2 z}{\partial x^2}(x,t) \\
\frac{\partial z}{\partial x} (0,t)=0, &&\frac{\partial z}{\partial x} (1,t)=0
\eean
so the open loop eigenvalues are $\lambda_n=-n^2\pi^2$ for $n=0,1,2,\ldots$
and the orthonormal eigenvectors are
\bea \label{olef}
\phi^0(x)=1,&&\phi^n(x)=\sqrt{2} \cos n\pi x
\eea
for  $n=1,2,\ldots$.	   Notice that $\lambda_0=0$ so the open loop system
 is only neutrally stable.  The rest of the eigenvalues are rapidly going to $-\infty$.

 We wish to stabilize the rod to some uniform temperature which we conveniently take to be $z=0$ by using a Linear Quadartic Regulator (LQR).
 We choose some $Q(x)\ge 0$  and a $m\times m$ positive definite matrix $R>0$ 
 and we seek to minimize
 \bea \label{crit}
\int_0^\infty \int_0^1 z(x,t)Q(x)z(x,t)\ dx+ u'(t)Ru(t)\ dt
 \eea
 subject to the above dynamics.

 Let $P(x_1,x_2)$ be any symmetric  function, $P(x_1,x_2)=P(x_2,x_1)$, which is continuous on the unit square $S=[0,1]^2$ and suppose there is a control trajectory $u(t)$ such that $z(x,t)\to 0$ as $t\to \infty $.  Then
 by the Fundamental Theorem of Calculus 
 \bea \nonumber
 &&0=\iint_{S}z^0(x_1)P(x_1,x_2)z^0(x_2)\ dA +\int_0^\infty \iint_{S} {d\over dt}\left(z(x_1,t)P(x_1,x_2)z(x_2,t)\right)\ dA \ dt\\
 &&0=\iint_{S}z^0(x_1)P(x_1,x_2)z^0(x_2)\ dA+\int_0^\infty \iint_{S} \frac{\partial^2 z}{\partial x_1^2}(x_1,t)P(x_1,x_2)z(x_2,t)\ dA \ dt\nonumber \\
 && +\int_0^\infty \iint_{S} z(x_1,t)P(x_1,x_2)\frac{\partial^2 z}{\partial x_2^2}(x_2,t)\ dA \ dt \label{FTC}
 \eea
 where $dA=dx_1dx_2$.
 
 We assume that $P(x_1,x_2)$ satisfies Neumann boundary conditions in each variable
 \bean
 \frac{\partial P}{\partial x_1}(0,x_2)&=\ 0 \ =&  \frac{\partial P}{\partial x_1}(1,x_2)\\
  \frac{\partial P}{\partial x_2}(x_1,0)&=\ 0 \ =&  \frac{\partial P}{\partial x_2}(x_1,0)
 \eean

 Now we formally integrate by parts twice with respect to $x_1$ on each subinterval $[\xi_k,\xi_{k+1}]$ ignoring the fact that we only assumed 
 $P(x_1,x_2)$ is continuous,
 \bea
 &&\int_{\xi_k}^{\xi_{k+1}} \frac{\partial^2 z}{\partial x_1^2}(x_1,t)P(x_1,x_2)z(x_2,t)\ dx_1=\int_{\xi_k}^{\xi_{k+1}} z(x_1,t)\frac{\partial^2 P}{\partial x_1^2}(x_1,x_2) z(x_2,t)\ dx_1\nonumber\\
 && +\left[\frac{\partial z}{\partial x_1}(x_1,t)P(x_1,x_2) z(x_2,t)\right]_{x_1=\xi_k^+}^{x_1=\xi_{k+1}^-}  -\left[z(x_1,t)\frac{\partial P}{\partial x_1}(x_1,x_2) z(x_2,t)\right]_{x_1=\xi_k^+}^{x_1=\xi_{k+1}^-} \nonumber\\
  \label{int_k}
\eea

Since we assumed that $\frac{\partial z}{\partial x_1}(x,t)=0$ off of $[0,1]$ we have
\bean
\frac{\partial z}{\partial x_1}(\xi_1^-,t)&=& \frac{\partial z}{\partial x_1}(\xi_m^+,t)\ =\ 0
\eean 
We sum (\ref{int_k})  over $k=1,\ldots,m-1$ and obtain \bean
&& \int_0^1 \frac{\partial^2 z}{\partial x_1^2}(x_1,t)P(x_1,x_2)z(x_2,t)\ dx_1\\
&&=\sum_{k=1}^{m-1} \int_{\xi_k}^{\xi_{k+1}} \frac{\partial^2 z}{\partial x_1^2}(x_1,t)P(x_1,x_2)z(x_2,t)\ dx_1\\
&&=\int_0^1 z(x_1,t)\frac{\partial^2 P}{\partial x_1^2}(x_1,x_2) z(x_2,t)\ dx_1+\sum_{k=1}^{m}\left( \frac{\partial z}{\partial x_1}(\xi^-_{k},t)-
\frac{\partial z}{\partial x_1}(\xi^+_{k},t)\right)P(\xi_{k},x_2) z(x_2,t)\\
&&=\int_0^1 z(x_1,t)\frac{\partial^2 P}{\partial x_1^2}(x_1,x_2) z(x_2,t)\ dx_1+\sum_{k=1}^{m}\beta_k u_k(t)P(\xi_{k},x_2) z(x_2,t)
\eean

Similarly
\bean
&& \int_0^1 z(x_1,t)P(x_1,x_2)\frac{\partial^2 z}{\partial x_2^2}(x_2,t)\ dx_2\\
&&=\int_0^1 z(x_1,t)\frac{\partial^2 P}{\partial x_2^2}(x_1,x_2) z(x_2,t)\ dx_2+\sum_{k=1}^{m}z(x_1,t)P(x_1,\xi_{k}) \beta_k  u_k(t)
\eean

We plug these into (\ref{FTC}) and obtain the identity
\bea \label{FTC1}
&&0=\iint_{S}z^0(x_1)P(x_1,x_2)z^0(x_2)\ dA\\
&&+\int_0^\infty \iint_{S} z(x_1,t)\nabla^2  P(x_1,x_2)z(x_2,t)\nonumber \\
&&+\sum_{k=1}^{m}\beta_k u_k(t)P(\xi_{k},x_2) z(x_2,t) +\sum_{k=1}^{m}z(x_1,t)P(x_1,\xi_{k}) \beta_k  u_k(t)\ dA \ dt \nonumber
 \eea
where $\nabla^2  P(x_1,x_2)$ denotes the two dimensional Laplacian of $P$.

We add the right side of this identity (\ref{FTC1}) to the criterion (\ref{crit}) to be minimized to get the equivalent criterion
\bea \label{crit1}
&&\iint_S z^0(x_1)P(x_1,x_2)z^0(x_2)\ dA  \\
&&+\int_0^\infty \iint_{S}z(x_1,t)\delta(x_1-x_2)Q(x_1)z(x_2,t)+ u'(t)Ru(t) \nonumber\\
&&+ z(x_1,t)\nabla^2  P(x_1,x_2)z(x_2,t)\nonumber \\
&&+\sum_{k=1}^{m}\beta_k u_k(t)P(\xi_{k},x_2) z(x_2,t) +\sum_{k=1}^{m}z(x_1,t)P(x_1,\xi_{k})   \beta_k u_k(t)\ dA \ dt \nonumber
 \eea
 where $\delta(x)$ is the Dirac delta function.

We would like to find a function $K(x)$ taking values in $\Real^{m\times 1}$ such that 
\bea
&&\iint_S \left(u(t)-K(x_1)z(x_1,t)\right)'R\left(u(t)-K(x_2)z(x_2,t)\right) dA
\\
&&=\iint_Sz(x_1,t)\delta(x_1-x_2)Q(x_1)z(x_2,t)+ u'(t)Ru(t) \nonumber+ z(x_1,t)\nabla^2  P(x_1,x_2)z(x_2,t)\nonumber \\
&&+\sum_{k=1}^{m}\beta_k u_k(t)P(\xi_{k},x_2) z(x_2,t) +\sum_{k=1}^{m}z(x_1,t)P(x_1,\xi_{k}) \beta_k  u_k(t) \ dA\nonumber
\eea

Clearly the terms quadratic in $u(t)$ agree so we we compare terms bilinear in $ z(x_2,t)$
and $u_k(t)$.  This yields 
\bea  \label{k_j}
-\sum_{j=1}^m R_{k,j}K_j(x_2) =\beta_kP(\xi_k,x_2)
\eea
Let $\beta$ be the $m\times m$ diagonal matrix with diagonal entries $\beta_1,\ldots,\beta_m$ and 
$\bar{P}(x_2)=[P(\xi_1,x_2),\ldots,P(\xi_m,x_2)]'$
then (\ref{k_j}) becomes the equation 
\bean
-RK(x_2) &=& \beta \bar{P}(x_2)
\eean
Therefore we define 
\bean
K(x) &=& -R^{-1}\beta \bar{P}(x)
\eean

Next we compare terms bilinear in $ z(x_1,t)$ and $ z(x_2,t)$,
\bean
&&\iint_S z(x_1,t)\left(\nabla^2  P(x_1,x_2)+\delta(x_1-x_2)Q(x_1)\right)z(x_2,t)\ dA\\
&&=
\iint_S \left(\beta \bar{P}(x_1)z(x_1,t)\right)'R^{-1}\beta \bar{P}(x_2)z(x_2,t)\ dA
\eean
This will hold if $P(x_1,x_2)$ is a solution  what we call a Riccati  PDE,
\bea \label{rpde}
\nabla^2  P(x_1,x_2)+\delta(x_1-x_2)Q(x_1)= \left(\beta \bar{P}(x_1)\right)'R^{-1}\beta \bar{P}(x_2)
\eea

Recall we only assumed that $P(x_1,x_2)$  is continuous on the unit square so the Riccati
PDE must be interpreted in the weak sense.  This means that given any $C^2$ function $\psi(x)$ satisfying 
Neumann boundary conditions it must be true that
\bean
0=\iint_S \left(\nabla^2  P(x_1,x_2)+\delta(x_1-x_2)Q(x_1)- \left(\beta \bar{P}(x_1)\right)'R^{-1}\beta \bar{P}(x_2)\right)
\psi(x_1)\psi(x_2)\ dA
\eean 
where the integral is evaluated using integration by parts.

We assume $Q(x)=q>0 $ then formally
\bean
Q(x)\delta(x_1-x_2)&=& q \sum_{n_1,n_2=0}^\infty \delta_{n_1,n_2} \cos n_1\pi x_1\cos n_2\pi x_2\
\eean
where $\delta_{n_1,n_2} $ is the Kronecker delta.

We assume that  $P(x_1,x_2)$ has a similar expansion 
\bean
P(x_1,x_2)&=&  \sum_{n_1,n_2=0}^\infty \delta_{n_1,n_2} P^{n_1,n_2}\cos n_1\pi x_1\cos n_2\pi x_2
\eean
and we define 
\bean
\gamma^{n,n}&=&\bmt \cos n\pi  \xi_1\\ \vdots\\ \cos n\pi  \xi_m\emt' \beta' R^{-1} \beta \bmt \cos n\pi \xi_1\\ \vdots\\ \cos n\pi \xi_m\emt
 \eean
Then the Riccati  PDE (\ref{rpde}) reduces to a sequence of quadratic equations
\bean
-2n^2\pi^2 P^{n,n}+q &=& \gamma^{n,n}\left(P^{n,n}\right)^2
\eean
with  roots
\bea \label{Pnn}
P^{n,n}&=& -n^2\pi^2+\sqrt{n^4\pi^4+  \gamma^{n,n} q}
\eea

The roots are positive since $q>0$ and $\gamma^{n,n}> 0 $ if at least one $\beta_i\ne 0$.  But as we now show they are going to zero like ${1\over n^2}$.  The Mean Value Theorem implies
there exists an $s$ between $n^4\pi^4 $ and $n^4\pi^4+  \gamma^{n,n} q$
such that 
\bean
P^{n,n}&=&{1\over 2\sqrt{s}}  \gamma^{n,n} q
\eean
The maximum of ${1\over 2\sqrt{s}} $  between $n^4\pi^4 $ and $n^4\pi^4+  \gamma^{n,n} q$
occurs at $s=n^4\pi^4 $ so we get the estimate
\bean
P^{n,n}&\le &{1\over 2n^2\pi^2}  \gamma^{n,n} q
\eean 
for $n>0$. 
Hence the series
\bea \label{Psum}
P(x_1,x_2) &=& \sum_{n=0}^\infty P^{n,n}\cos n\pi x_1\cos n\pi x_2 
\eea
converges uniformly to a continuous function which is a weak solution of the Riccati PDE (\ref{rpde}).

\section{Example One}
We first consider a simple example with heating/cooling at both ends of the rod, $m=2$, $\xi_1=0$,
$\xi_2=1$, $\beta_1=\beta_2=1$, $q=1$, $R=[1,0;0,1]$ and $\gamma^{n,n}=2$.
Then
\bean
K_1(x)&=& - \sum_{n=0}^\infty P^{n,n}\cos n\pi x \\
K_2(x)&=&- \sum_{n=0}^\infty (-1)^n P^{n,n}\cos n\pi x 
\eean
The closed loop "boundary" conditions are
\bean
\frac{\partial z}{\partial x}(0^+,t)&=& \sum_{n=0}^\infty\int_0^1 P^{n,n}\cos n\pi x \ z(x,t)\ dx \\
\frac{\partial z}{\partial x}(1^-,t)&=& -\sum_{n=0}^\infty\int_0^1 (-1)^n P^{n,n}\cos n\pi x \ z(x,t)\ dx
\eean
Notice that these nonstandard "boundary" conditions are linear. A straightforward calculation
from (\ref{Pnn})  yields
$P^{0,0}=1.4142$,   $P^{1,1}=0.1008$, $P^{2,2}=    0.0253$,  $P^{3,3}=    0.0113$, $P^{4,4}=    0.0063$ and $P^{5,5}=    0.0041$.

An eigenfunction of the Laplacian under any linear   "boundary" conditions
is of the form
\bean
\psi(x)&=& a\cos \nu  x+b\sin \nu x 
\eean 
for some $\nu$ and the corresponding eigenvalue is $-\nu^2$. 
 
For this example the "boundary" conditions  are
\bea \label{BC1}
b \nu &=&\sum_{n=0}^\infty \int_0^1 P^{n,n} \cos n\pi x \left( a \cos \nu x+b\sin \nu x\right)\ dx \\ \nonumber
\\
  a \nu \sin \nu - b \nu \cos \nu &=& \sum_{n=0}^\infty\int_0^1(-1)^n P^{n,n} \cos n\pi x \left( a \cos \nu x+b\sin \nu x\right)\ dx\nonumber \\
  \label{BC2}
\eea

Due to the symmetry of the problem the eigevectors must be invariant under the action of replacing $x$ by $1-x$.
  So the eigenvectors are either $\psi(x)=\cos 2k\pi x$ or $\psi(x)=\sin (2k+1)\pi x$ for some nonegative integer $k$.
  
  If $k=0$ and $\psi(x)=1$  the "boundary" conditions  become
\bean
0&=&\sum_{n=0}^\infty \int_0^1 P^{n,n} \cos n\pi x \ dx\ =\ 1\\
\\
0&=&\sum_{n=0}^\infty \int_0^1 (-1)^n P^{n,n} \cos n\pi x \ dx \ =\ 1
\eean
so $\psi(x)=1$ is not an eigenfunction.

If $k>0$ and $\psi(x)=\cos 2k\pi x$ the "boundary" conditions  become
\bean
2k\pi &=&\sum_{n=0}^\infty \int_0^1 P^{n,n} \cos n\pi x \cos 2k\pi x\ dx\ =\ {1\over 2}\\ 
\\
-2k\pi &=&\sum_{n=0}^\infty \int_0^1 (-1)^n P^{n,n} \cos n\pi x \cos 2k\pi x \ dx \ =\ {1\over 2}
\eean
so $\psi(x)=\cos 2k\pi x$ is not an eigenfunction.

If $k\ge 0$ and $\psi(x)=\sin (2k+1) \pi x$ the "boundary" conditions  are
\bean
(2k+1)\pi &=& \sum_{n=0}^\infty \int_0^1 P^{n,n} \cos n\pi x \sin (2k+1) x\ dx  \\
\\
(2k+1)\pi &=& \sum_{n=0}^\infty \int_0^1 (-1)^n P^{n,n} \cos n\pi x \sin (2k+1) x\ dx 
\eean
If $n$ is odd then $n\pm (2k+1)$ are even
so
\bean
&&\int_0^1   \cos n\pi x \sin (2k+1) \pi x\ dx \\
&&={1\over 2} \int_0^1 \sin (n+ (2k+1))\pi x +\sin (n-(2k+1))\pi x \ dx=0
\eean
This shows that the two boundary conditions are identical.
So  the  closed loop eigenfunctions are  $\psi_k(x)=\sin (2k+1) \pi x$ and the closed loop eigenvalues are $\mu_k =- (2k+1)^2 \pi^2$
 for $k=0,1,2,\ldots$.

 Being able to heat/cool the rod at both ends is a big improvement over being able to heat/cool the rod at just one end.
As we saw in \cite{Kr20b} for control only at one end the least stable closed loop eigenvalue is $-1.0409$.   For control at both ends the least stable closed loop eigenvalue is $-\pi^2=-9.8696$.

\section{Example Two}
We assume that the rod can be heated/cooled at both endpoints and also at the midpoint,
$m=3$, $\xi_1=0, \ \xi_2=0.5,\ \xi_3=1$, $\beta_1=1,\ \beta_2=2,\ \beta_3=1$, $q=1$, $R=[1,0,0;0,1,0;0,0,1]$  and
$\gamma^{n,n}=6$ if $n$ is even and $\gamma^{n,n}=2$ if $n$ is odd.

The $P^{n,n}$ are given by (\ref{Pnn}).
The optimal feedback gains are
\bea
\label{K1}
K_1(x)&=& -\sum_{n=0}^\infty P^{n,n} \cos n\pi x\\
\label{K2}
K_2(x)&=&-2\sum_{k=0}^\infty  (-1)^{k} P^{2k,2k}\cos 2k\pi x\\
\label{K3}
K_3(x)&=& -\sum_{n=0}^\infty(-1)^n P^{n,n} \cos n\pi x
\eea

We assume that a closed loop eigenfunction $\psi(x)$ is composed of different sinusoids on $[0,0.5]$ and $[0.5,1]$
with a common frequency $\nu$.  Because the system is
symmetric with respect to replacing $x$ with $1-x$ we expect a closed eigenfunction to reflect this symmetry,
\bean
\psi (x)&=& \left\{\begin{array}{ccc} a \cos \nu  x+b  \sin \nu x& \mbox{ if } & 0\le x\le 0.5\\ \\
a \cos \nu (1-x)+b  \sin \nu (1-x)& \mbox{ if } &0.5\le x\le 1
\end{array}\right.
\eean

Notice such a solution immediately satisfies the continuity condition,
\bea \label{cont}
\psi (0.5^-)&=&\psi (0.5^+)
\eea
and
\bea \label{jump}
\frac{\partial \psi}{\partial x}(0.5^-)&=&-\frac{\partial \psi}{\partial x}(0.5^+)
\eea

The first closed loop  "boundary" condition is
\bean
b &=& \sum_{n=0}^\infty \int_0^{0.5} P^{n,n} \cos n\pi x\left(a \cos \nu  x+b  \sin \nu x\right)\ dx
\\
&& +\sum_{n=0}^\infty \int_{0.5}^1 P^{n,n} \cos n\pi x\left(a \cos \nu (1-x)+b\sin (1-x)\right)\  dx\\
\eean

Now $\cos n\pi (1-x)=(-1)^n \cos n\pi x$ so 
\bean
\int_{0.5}^1 P^{n,n} \cos n\pi x\left(a \cos \nu (1-x)+b\sin \nu (1-x)\right)\  dx
\\
=\int_0^{0.5} (-1)^n P^{n,n} \cos n\pi x\left(a \cos \nu x+b\sin \nu x\right)\  dx
\eean
Hence the first closed loop  "boundary" condition becomes
\bean
b&=& 2\sum_{k=0}^\infty  \int_0^{0.5} P^{2k,2k} \cos 2k\pi x\left(a \cos \nu  x+b  \sin \nu x\right)\ dx
\eean

The second closed loop  "boundary" condition is
\bean
\frac{\partial \psi}{\partial x}(0.5^-) -\frac{\partial \psi}{\partial x}(0.5^+) &=& \int_0^1 K_2(x)
\psi(x)\ dx
\eean
From (\ref{jump}) and (\ref{K2})
we obtain
\bean 
&&2\nu\left (-a \sin 0.5  \nu +b\cos  0.5  \nu\right)\\
&&= -2\sum_{k=0}^\infty  \int_0^{0.5}(-1)^{k} P^{2k,2k}\cos 2k\pi x\left( a\cos \nu x+b\sin \nu x\right)\ dx\\
&&-2\sum_{k=0}^\infty  \int_{0.5}^1(-1)^{k} P^{2k,2k}\cos 2k\pi x \left( a\cos \nu (1-x)+b\sin \nu (1-x)\right)\ dx\\
&&=-4\sum_{k=0}^\infty  \int_0^{0.5}P^{2k,2k}\cos 2k\pi x\left( a\cos \nu x+b\sin \nu x\right)\ dx
\eean

If $0.5  \nu$ is an odd integer, i.~e.~, $\nu=4j+2$, then the two boundary conditions are identical so the closed loop
eigenvalues and eigenvectors are $\mu_j=-(4j+2)^2\pi^2$ and $\psi_j(x)=\sin (4j+1)\pi x$ for $j=0,1,2,\ldots$ and $x\in [0,0.5]$.
In particular the least stable closed loop eigenvalue is $-4\pi^2=-39.4784$ so being able to heat/cool in the middle has a big impact.

\section{Kalman Filtering of the Heat Equation with Point Observations}
In the previous sections we constructed optimal feedbacks to control the heat equation.
These feedbacks assumed that the full state $z(x,t)$ is known at every $ x\in [0,1]$ and $t\ge 0$.  But in practice we may only be able to
measure the temperature at a finite number of points $z(\zeta_1,s), z(\zeta_2,s), z(\zeta_p,s)$
where $ 0\le \zeta_1<\zeta_2<\ldots \zeta_p\le 1$ for $-\infty< s\le t$ and these measurements may be corrupted by noise.

Our model for the measured but uncontrolled rod is
\bean
\frac{\partial z}{\partial t}(x,t)&=& \frac{\partial^2 z}{\partial x^2}(x,t) +B v(t)\\
y_i(t)&=&C_i z(\zeta_i,t)+  D_i w_i(t)\\
\frac{\partial z}{\partial x}(0,t)=0,&&\frac{\partial z}{\partial x}(1,t)=0
\eean
where $v(t),\ w_i(t)$ are independent white Gaussian noise processes.  Without loss of generality we can assume $C_i>0$ for $i=1,\ldots,p$.

  We
wish to construct an estimate $\hat{z}(x,t)$ for all $x\in[0,1]$ based on the past measurements, $y_i(s),\ s\le t$ .  We assume the estimates 
are linear functionals of  the past  observations of the form 
\bean
\hat{z}(x,t)&=& \int_{-\infty}^t \sum_{i=1}^p {\cal K}_i(x,s-t)y_i(s) \ ds
\eean
 Since we are taking
measurements for $-\infty<s\le t$ we expect the filter to be stationary.  Therefore it suffices to solve the problem for $t=0$ and we only need to consider $y_i(s)$ for $-\infty<s\le 0$.

Given a possible set of filter gains ${\cal K}_i(x,s)$  we define ${\cal H}x,x_1,s)$ as the solution of a driven backward generalized heat equation
\bean
 \frac{\partial {\cal H}}{\partial s}(x,x_1,s)\ &=& -\frac{\partial^2 {\cal H}}{\partial x_1^2}(x,x_1,s) +\sum_{i=1}^p {\cal K}_i(x,s)C_i \delta(x_1-\zeta_i) \\
{\cal H}(x,x_1,0)&=&\delta(x-x_1)
\eean
where ${\cal H}(x,x_1,s)$ satisfies Neumann boundary conditions with respect to  $x$ and $x_1$.

Then
\bean
\hat{z}(x,0)&=& \int_{-\infty}^0 \sum_{i=1}^p \int_0^1 {\cal K}_i(x,s) C_i \delta(x_1-\zeta_i) z(x_1,s)\ dx_1\ ds\\&&+\int_{-\infty}^0 \sum_{i=1}^p {\cal K}_i(x,s) D_iw_i(s)\ ds\\
&=&\int_{-\infty}^0 \int_0^1 \left(\frac{\partial {\cal H}}{\partial s}(x,x_1,s)+ \frac{\partial^2{\cal H}}{\partial x_1^2}(x,x_1,s)\right)z(x_1,s) \ dx_1 \\
&&+\sum_{i=1}^p {\cal K}_i(x,s)D_iw_i(s)\ ds
\eean
We integrate by parts with respect to $s$ assuming ${\cal H}(x,x_1,s)\to 0$ as $ s\to -\infty$ and obtain
\bean
\hat{z}(x,0)-z(x,0)&=&
-\int_{-\infty}^0 \int_0^1  {\cal H}(x,x_1,s) \frac{\partial z}{\partial s}(x_1,s)\\
&&-\frac{\partial^2 {\cal H}}{\partial x_1^2}(x,x_1,s)z(x_1,s) \ dx_1\ ds \\
&&+\int_{-\infty}^0\sum_{i=1}^p {\cal K}_i(x,s) D_iw_i(s)\ ds
\eean
\bean
\hat{z}(x,0)-z(x,0)&=&
 -\int_{-\infty}^0 \int_0^1  {\cal H}(x,x_1,s)\left( \frac{\partial^2 z}{\partial x_1^2}(x_1,s)+ Bv(s)\right)\\
&& -\frac{\partial^2 {\cal H}}{\partial x_1^2}(x,x_1,s)z(x_1,s) \ dx_1\ ds \\
&&+\int_{-\infty}^0\sum_{i=1}^p {\cal K}_i(x,s)D_iw_i(s)\ ds
\eean

We  integrate the second term on the right by parts twice with respect to $x_1$ using the Neumann boundary conditions and obtain 
\bean
\hat{z}(x,0)-z(x,0)&=&
 -\int_{-\infty}^0  {\cal H}(x,x_1,s)Bv(s) \ dx_1\ ds \\
&&+\int_{-\infty}^0\sum_{i=1}^p {\cal K}_i(x,s) D_iw_i(s)\ ds
\eean
 Since $v(s)$ and $w_i(s)$ are independent white Gaussian noise processes, the estimation error $\tilde{z}(x,0)=z(x,0)-\hat{z}(x,0)$ has error variance
 \bea \label{EE} \nonumber
{\rm E}\left(\tilde{z}(x,0)\left(\tilde{z}(x,0)\right)'\right)&=& \int_{-\infty}^0 \int_0^1  {\cal H}(x,x_1,s) B^2{\cal H}(x,x_1,s) \ dx_1 \ dt\\
&&+ \sum_{i=1}^p \int_{-\infty}^0 {\cal K}_i(x,s) D^2_i {\cal K}_i(x,s) \ ds
 \eea

For each $x\in [0,1] $  this is a backward  LQR optimal control problem with infinite dimensional state $x_1\to {\cal H}(x,x_1,s)$ and $m$ dimensional  control ${\cal K}_i(x,s), \ i=1,\ldots,p$. 
Our goal is to minimize for each $x$ the error variance  subject to the backward dynamics 
\bean
 \frac{\partial {\cal H}}{\partial s}(x,x_1,s)\ &=& -\frac{\partial^2 {\cal H}}{\partial x_1^2}(x,x_1,s) +\sum_{i=1}^p {\cal K}_i(x,s)C_i \delta(x_1-\zeta_i) 
\eean
for $-\infty<s\le 0$.   The terminal condition is
\bean
{\cal H}(x,x_1,0)&=&\delta(x-x_1)
\eean
and the  boundary conditions are 
\bean
\frac{\partial {\cal H}}{\partial x_1}(x,0,s)=0,&& \frac{\partial {\cal H}}{\partial x_1}(x,1,s)=0
\eean

For any $x\in[0,1]$, let $P(x,x_1,x_2)$ be any continuous function symmetric with respect to $x_1,x_2$,
$P(x,x_1,x_2)=P(x,x_2,x_1)$, then given that $  {\cal H}(x,x_1,s)\to 0$ as $s\to -\infty$ then
\bean
&&0=\iint_S  P(x,x_1,x_2){\cal H}(x,x_1,0){\cal H}(x,x_2,0) \ dA \\
&&-\int_{-\infty}^0 \iint_S {d\over ds} P(x,x_1,x_2){\cal H}(x,x_1,s){\cal H}(x,x_2,s)\ dA \ ds \\
&&=- \int_{-\infty}^0 \iint_S P(x,x_1,x_2) {\partial {\cal H}\over \partial s}(x,x_1,s){\cal H}(x,x_2,s)\ dA \ ds 
 \\
&& -\int_{-\infty}^0 \iint_S  P(x,x_1,x_2){\cal H}(x,x_1,s){\partial {\cal H}\over \partial s}(x,x_2,s)\ dA \ ds 
\eean
 where $dA=dx_1dx_2$ and $S$ is the unit square in $x_1,x_2$ plane.
 
 We plug in the dynamics of ${\cal H}$ to get
\bean
&&0=\iint_S  P(x,x_1,x_2){\cal H}(x,x_1,0){\cal H}(x,x_2,0) \ dA \\
&&- \int_{-\infty}^0 \iint_S P(x,x_1,x_2)\left(-\frac{\partial^2 {\cal H}}{\partial x_1^2}(x,x_1,s)+\sum_{i=1}^p {\cal K}_i(x,s)C_i\delta(x_1-\zeta_i)\right)\\
&&\times {\cal H}(x,x_2,s)\ dA \ ds\\
&&-\int_{-\infty}^0 \iint_S P(x,x_1,x_2){\cal H}(x,x_1,s)\\
&& \times \left(-\frac{\partial^2 {\cal H}}{\partial x_2^2}(x,x_2,t)+ \sum_{i=1}^p {\cal K}_i(x,s)C_i\delta(x_2-\zeta_i) \right) \ dA \ ds\
\eean
then we integrate by parts twice to get
\bean
&&0=\iint_S  P(x,x_1,x_2,0){\cal H}(x,x_1,s){\cal H}(x,x_2,s) \ dA \\
&&-\int_{-\infty}^0 \iint_S -\nabla^2 P(x,x_1,x_2,s){\cal H}(x,x_1,s){\cal H}(x,x_2,s)\ dA \ ds
\\
&&-\int_{-\infty}^0 \int_0^1 {\cal H}(x,x_2,s)\sum_{i=1}^p P(x,\zeta_i,x_2){\cal K}_i(x,s)C_i\ dx_2 \ ds
\\
&&-\int_{-\infty}^0 \int_0^1 {\cal H}(x,x_1,s) \sum_{i=1}^p P(x,x_1\zeta_i){\cal K}_i(x,s)C_i\ dx_1 \ ds
\eean
where $\nabla^2$ is the two dimensional Lagrangian with respect to $x_1$ and $x_2$.

We add the right side of this identity to the estimation error variance (\ref{EE}) to get an equivalent
quantity to be minimized
\bean
&& \int_{-\infty}^0 \iint_S \sum_{j=1}^m {\cal H}(x,x_1,s) B^2 \delta(x_1-x_2){\cal H}(x,x_2,s) \ dA \ dt\\
&&+ \sum_{i=1}^p \int_0^\infty {\cal K}_i(x,s) D^2_i {\cal K}_i(x,s)\ ds \\
 &&-\iint_S  P(x,x_1,x_2){\cal H}(x,x_1,s){\cal H}(x,x_2,s) \ dA \\
&&- \int_{-\infty}^0 \iint_S- \nabla^2 P(x,x_1,x_2){\cal H}(x,x_1,s){\cal H}(x,x_2,s)\ dA \ ds
\\
&&-\int_{-\infty}^0 \int_0^1{\cal H}(x,x_2,s) \sum_{i=1}^p P(x,\zeta_i,x_2){\cal K}_i(x,s)C_i\ dx_2 \ ds
\\
&&-\int_{-\infty}^0 \int_0^1 {\cal H}(x,x_1,s)\sum_{i=1}^p P(x,x_1,\zeta_i){\cal K}_i(x,s)C_i\ dx_1 \ ds
\eean
 
 For each $i=1,\ldots, p$ and each $x\in[0,1]$ we would to choose $L_i(x,x_1,s)$ so the time integrand of the quantity 
 to be minimized is a perfect square of the form
 \bean
\int_S \sum_{i=1}^p \left({\cal K}_i(x,s)-L_i(x,x_1){\cal H}(x,x_1,s)\right) D^2_i \left({\cal K}_i(x,s)-L_i(x,x_2){\cal H}(x,x_2,s)\right)\ dA\ dt
\eean

Clearly the terms quadratic in ${\cal K}_i(x,s)$ match up so we compare terms bilinear in ${\cal K}_i(x,s)$ and ${\cal H}(x,x_2,s)$,
\bean
&&\int_0^1\sum_{i=1}^p {\cal K}_i(x,s)D^2_i L_i(x,x_2){\cal H}(x,x_2,s)\ dx_2\\
&&= \int_0^1 \sum_{i=1}^p {\cal H}(x,x_2,s)P(x,\zeta_i,x_2){\cal K}_i(x,s)C_i\ dx_2
\eean
This will hold if 
\bean
D^2_i L_i(x,x_2){\cal H}(x,x_2,s)&=& {\cal H}(x,x_2,s)P(x,\zeta_i,x_2)C_i
\eean
for $i=1,\ldots,p$ so we define 
\bean
L_i(x,x_1)&=&D_i^{-2} P(x,x_1,\zeta_i)C_i
\eean

   Then we compare terms bilinear  in ${\cal H}(x,x_1,s)$ and ${\cal H}(x,x_2,s)$ and obtain
   \bean
   &&\iint_S  {\cal H}(x,x_1,s)B^2\delta(x_1-x_2){\cal H}(x,x_2,s) \\
&&-  \iint_S -\nabla^2 P(x,x_1,x_2){\cal H}(x,x_1,s){\cal H}(x,x_2,s)\ dA
\\
&&=\iint_S  \sum_{i=1}^p L_i(x,x_1){\cal H}(x,x_1,s)D^2L_i(x,x_2){\cal H}(x,x_2,s)\ dA\\
&&=\iint_S  \sum_{i=1}^p {\cal H}(x,x_1,s)P(x,\zeta_i,x_1)C^2_iD^{-2} P(x,\zeta_i,x_2)\ dA 
\eean

So we are looking for a weak solution to what we call the filter Riccati PDE,
\bean
0&=&- \nabla^2 P(x,x_1,x_2)+ B^2\delta(x_1-x_2)
\\&&+\sum_{i=1}^p P(x,x_1,\zeta_i)C^2_i\ D^{-2} P(x,\zeta_i,x_2)
\eean

We guess that $P(x,x_1,x_2)$ has an  expansion 
\bean
P(x,x_1,x_2)&=& \sum_{n=0}^\infty P^{n,n}(x) \phi_n(x_1)\phi_n( x_2)
\eean
where $ \phi_n(x)$ are the orthonomal eigenfunctions of the Laplacian
under Neumann boundary conditions (\ref{olef}).
We plug this into filter Riccati PDE and we obtain for $n=0,1,2,\ldots$ the equations
\bean
\left(2n^2 \pi^2\right) P^{n,n}(x) +B^2 \delta_{0,n} +\sum_{i=1}^p C_i^2 \ D_i^{-2}\left(P^{n,n}(x)\right)^2
\phi_n( \zeta_i)
\eean
One solution to these equations   is $P^{n,n}(x)=0$ if  $n>0$  and
\bean
P^{0,0}(x)&=& P^{0,0}\ =\ \sqrt{B^2\over \sum_{i=1}^p C_i^2 \ D_i^{-2}}
\eean
So $P(x,x_1,x_2)=P^{0,0}$ and 
the fact that it does not depend on $x$ is not surprising as the coefficient $B$
of the driving noise is constant.  If the driving noise had some spatial variation
we suspect that $P(x,x_1,x_2)$ would vary with $x$.

Then $L_i(x,x_1) $ is a constant and 
\bean
L_i(x,x_1)&=&L_i\ = \  P^{0,0}C_i\\
{\cal K}_i(x,s) &=& \int_0^1 L_i{\cal H}(x,x_1,s) \ dx_1\ = \ P^{0,0}C_i \int_0^1 {\cal H}(x,x_1,s) \ dx_1
\eean 
Since $P^{0,0}$ and $ C_i$ are both positive so is $L_i$ for $i=1,\ldots, p$.

Now  ${\cal H}(x,x_1,s)$ satisfies the backward linear partial differential equation 
\bea \label{bpde}
\frac{\partial {\cal H}}{\partial s}(x,x_1,s)&=&- \frac{\partial {\cal H}^2}{\partial x_1^2}(x,x_1,s)+\sum_{i=1}^p  L_i {\cal H}(x,x_1,s)
\eea 
subject to the terminal condition  ${\cal H}(x,x_1,0)=\delta(x-x_1)$ and Neumann boundary  conditions in both $x$
and $x_1$.
  
We assume that the solution to this PDE takes the form
\bean
{\cal H}(x,x_1,s)&=&\sum_{m,n=0}^\infty \gamma_{m,n}(s) \phi_m(x) \phi_n(x_1)
\eean
where $\phi_n(x)$ are the orthonormal open loop eigenfunctions (\ref{olef}).
 
We plug this into (\ref{bpde}) and we get a sequence of ODEs,
\bean
{ d\over ds} \gamma_{m,n}(s) &=&\left( n^2\pi^2+ \sum_{i=1}^p  L_i\right) \gamma_{m,n}(s)
\eean
The terminal  condition
$
{\cal H}(x,x_1,0)= \delta(x-x_1)
$
implies
\bean
\gamma_{m,n}(0)&=& \delta_{m,n}
\eean
so
\bean
\gamma_{m,n}(s)&=&\delta_{m,n}\exp\left(\left( n^2\pi^2+ \sum_{i=1}^p  L_i\right) s\right)
\eean
and 
\bea  \label{H}
{\cal H}(x,x_1,s)&=&\sum_{n=0}^\infty \gamma_{n,n}(s) \phi_n( x) \phi_n( x_1)
\eea
Notice that ${\cal H}(x,x_1,s)\to 0 $ as $s\to -\infty$.

Recall the estimate of the state at time zero is
\bean
\hat{z}(x,0)&=& \int_{-\infty}^0 \sum_{i=1}^p {\cal K}_i(x,s)y_i(s)\ ds
\eean
but this is a stationary estimator so the estimate of the state at time $t$ is
\bean
\hat{z}(x,t)&=& \int_{-\infty}^t \sum_{i=1}^p {\cal K}_i(x,s-t) y_i(s)\ ds\\
&=& \int_{-\infty}^t \sum_{i=1}^p L_i {\cal H}(x,x_1,s-t) y_i(s)\ ds
\eean
and
\bean
\frac{\partial \hat{z}}{\partial t} (x,t)&=&  \sum_{i=1}^p {\cal K}_i(x,0)y_i(t)+  \int_{-\infty}^t \sum_{i=1}^p \frac{\partial {\cal K}_i}{\partial t}(x,s-t) y_i(s)\ ds\\
&=&\sum_{i=1}^p {\cal K}_i(x,0)y_i(t)+  \int_{-\infty}^t \int_0^1\sum_{i=1}^p L_i \frac{\partial {\cal H}}{\partial t}(x,x_1,s-t) y_i(s)\ dx_1 \ ds\\
&=&\sum_{i=1}^p {\cal K}_i(x,0)y_i(t )+ \sum_{i=1}^p \int_{-\infty}^t \int_0^1 L_i \frac{\partial^2 {\cal H}}{\partial x^2_1} (x,x_1,s-t)y_i(s)\ dx_1 \ ds \\
&&- \sum_{i=1}^pL_i\int_{-\infty}^t \sum_{j=1}^p L_j C_j {\cal H}( x,\zeta_j,s) 
 y_i(s)\ ds
\eean

Notice
\bean
{\cal K}_i(x,0) &=& \int_0^1 L_i H(x,x_1,0) dx_1\ = \  \int_0^1 L_i \delta(x-x_1) dx_1\ = \ L_i
\eean

From (\ref{H}) we see that 
\bean
\frac{\partial^2 {\cal H}}{\partial x^2_1} (x,x_1,s-t)&=&\frac{\partial^2 {\cal H}}{\partial x^2} (x,x_1,s-t)
\eean
so 
\bean
\frac{\partial \hat{z}}{\partial t} (x,t)&=& \sum_{i=1}^p L_iy_i(t )+  \int_{-\infty}^t \int_0^1\sum_{i=1}^p L_i \frac{\partial^2 {\cal H}}{\partial x^2_1} (x,x_1,s-t)y_i(s)\ dx_1 \ ds \\
&&- \sum_{i=1}^pL_i\int_{-\infty}^t \sum_{j=1}^p L_j C_j {\cal H}( x,\zeta_j,s) 
 y_i(s)\ ds\\
 &=& \sum_{i=1}^p L_iy_i(t )+\frac{\partial^2 \hat{z}}{\partial x^2}(x,t)\\
 && -\sum_{j=1}^pL_j \int_{-\infty}^t \sum_{i=1}^p C_i{\cal H}(x,\zeta_j,s)y_i(s)\\
&=& \sum_{i=1}^p L_iy_i(t )+\frac{\partial^2 \hat{z}}{\partial x^2}(x,t)\\
 && -\sum_{j=1}^pL_j \int_{-\infty}^t \sum_{i=1}^p C_i{\cal H}(x,\zeta_j,s)y_i(s)\\
 &=& \sum_{i=1}^p L_iy_i(t )+\frac{\partial^2 \hat{z}}{\partial x^2}(x,t) -\sum_{j=1}^pL_j \hat{z}(t,\zeta_j)
 \eean

Now  we define
\bean
\hat{y}_i(t)&=& \hat{z}(t,\zeta_j)\\
\tilde{y}_i(t)&=&y_i(t)-\hat{y}_i(t)
\eean
The quantities $\tilde{y}_i(t)$ are called the innovations. 
So the dynamics of optimal estimator is a copy of original dynamics driven by the innovations
\bean
\frac{\partial \hat{z}}{\partial t} (x,t)&=&\frac{\partial^2 \hat{z}}{\partial x^2}(x,t) +\sum_{i=1}^pL_i \tilde{z}(t,\zeta_i)
\eean
The error  dynamics is  
\bea \label{errdyn}
\frac{\partial \tilde{z}}{\partial t} (x,t)&=&\frac{\partial^2 \tilde{z}}{\partial x^2}(x,t) -\sum_{i=1}^pL_i \tilde{z}(t,\zeta_i)
\eea
Because both $z(x,t)$ and $\hat{z}(x,t)$ satisfy Neumann boundary conditions so does $\tilde{z}(x,t)$.

From the form of the error dynamics we see that the closed loop eigenvalues $\eta_n$ and eigenfunctions $\theta_n(x)$ are weak solutions of the equations 
\bean
\frac{\partial^2 \theta_n}{\partial x^2}(x) -\sum_{i=1}^p L_i\delta (x-\zeta_i) \theta_n(x)&=&\eta_n \theta_n(x)
\eean
subject to Neumann boundary conditions.   
It is easy to see that that  $\theta_0(x)=1$ and $\eta_0=-\sum_i L_i<0$ since $L_i>0$.

For $n>0$ we expect that on each subinterval $\zeta_i< x< \zeta_{i+1}$ the eigenfunctions are a sinusoid 
of a given frequency $\tau_n  $. The corresponding eigenvalue is $-\tau_n^2 $.
  The eigenfuctions $\theta_n(x)$ are continuous at the mesurement points $\zeta_i$ but their derivatives jump
\bean
\frac{\partial \theta_n}{\partial x} \zeta_i(x^+)- \frac{\partial \theta_n}{\partial x} \zeta_i(x^-)&=& L_i \theta_n(\zeta_i)
\eean

Next we consider the point controlled and point measured system
\bean
\frac{\partial z}{\partial t}(x,t)&=& \frac{\partial^2 z}{\partial x^2}(x,t) +B v(t)\\
z(\xi_k^-,t)&=&z(\xi_k^+,t) \mbox{ for } k=2,\ldots, m-1\\
\beta_k u_k(t)&=&\frac{\partial z}{\partial x} (\xi_k^-,t) -\frac{\partial z}{\partial x} (\xi_k^+,t) \mbox{ for } k=1,\ldots, m
\\
y_i(t)&=&C_i z(\zeta_i,t)+  D_i w_i(t)\\
z(x,t) &=&  z^0(x)
\eean
We use the linear feedback on the state estimate to get the control inputs
\bean
\hat{u}(t)&=&\int_0^1 K(x) \hat{z}(x,t) \ dx 
\eean
where $\hat{z}(x,t)$ satisfies the Kalman filtering equation modified by the linear feedback on the state estimate,
\bean
\frac{\partial \hat{z}}{\partial t} (x,t)&=&\frac{\partial^2 \hat{z}}{\partial x^2}(x,t) +\sum_{i=1}^pL_i \tilde{z}(t,\zeta_i)\\
\hat{z}(\xi_k^-,t)&=&\hat{z}(\xi_k^+,t) \mbox{ for } k=2,\ldots, m-1\\
\beta_k \hat{u}_k(t)&=&\frac{\partial \hat{z}}{\partial x} (\xi_k^-,t) -\frac{\partial \hat{z}}{\partial x} (\xi_k^+,t) \mbox{ for } k=1,\ldots, m
\\
\hat{z}(x,t) &=&  \hat{z}^0(x)
\eean

By linearity the error  dynamics is  still 
\bean
\frac{\partial \tilde{z}}{\partial t} (x,t)&=&\frac{\partial^2 \tilde{z}}{\partial x^2}(x,t) -\sum_{i=1}^pL_i \tilde{z}(t,\zeta_i)\\
\tilde{z}(x,0)&=& \tilde{z}^0(x)\ =\ z^0(x)-\hat{z}^0(x)
\eean

We can consider the combined system in coordinates $z(x,t)$ and $\hat{z}(x,t)$ but it is more useful to consider it in coordinates $z(x,t)$ and $\tilde{z}(x,t)$ because in the latter coordinates the combined dynamics is block upper triangular, the error dynamics does not depend on $z(x,t)$.
In these coordinates the 
control is given by
\bean
\hat{u}(t)&=&\int_0^1 K(x) \left(z(x,t)- \tilde{z}(x,t)\right) \ dx 
\eean
so when $\tilde{z}(x,t)=0$ the 
dynamics of $z(x,t)$  takes the form of the original system under full state feedback.
This shows the spectrum of the LQG synthesis is the union of the spectrum of the original system under LQR full state feedback
and spectrum of the error dynamics of the Kalman filter.  So if these spectra are
in the open left half plane then the LQG synthesis is asymptotically stable, the rod goes to the desired temperature $z(x,t)\to 0$ and the estimation
error goes to zero $\tilde{z}(x,t)\to 0$

\section{Example 3}
We consider a Linear Quadratic Gaussian synthesis for the heat equation.  As in Example 2 we assume that there are three actuators, $m=3$, at $\xi_1=0,\ \xi_2=0.5, \  \xi_3=1$
with coefficients $\beta_1=1,\ \beta_2=2,\ \beta_3=1$ and the rest of the constants are  as in Example 2. 
We further assume that are  two sensors at $\zeta_1=0.25,\ \zeta_2=0.75$ with $C_1=C_2=1$.   The coefficient of the driving noise is $B=1$ and
the coefficients of measurement noise are $D_1=D_2=1$.  

Then
\bean
P^{0,0}&=& {\sqrt{2}\over 2}
\\
L_1&=& L_2\ = \  {\sqrt{2}\over 2}
\eean 

Then the zeroth order eigenfunction is $\theta_0(x)=1$ and the corresponding eigenvalue is $\eta_0=-\sqrt{2}$.  For $n>0$ the eigenfunctions are of the form
\bean
\theta_n(x)&=&\left\{ \begin{array}{ccccc} a_1\cos \tau_n  x  +b_1 \sin \tau_n x& 0\le x\le 0.25\\
a_2\cos \tau_n x  +b_2 \sin \tau_n  x& 0.25 \le x\le 0.75\\
a_3\cos \tau_n  x  +b_3 \sin \tau_n  x& 0.75\le x\le1
\end{array} \right.
\eean
By symmetry $\theta_n(x)=\theta_n(1-x)$ so only need to find $a_1,b_1,a_2,b_2$.  Since $\theta_n(x)$ satisfies Neumann boundary conditions
at $x=0$, $b_1=0$ and we can take $a_1=1$. 

On the interval $[0.25,0.75]$ the solution must be symmetric around $0.5$ so we make the change of coordinates $\bar{x} =0.5-x$ then
 it must be of the form
\bean
\theta_n(x)\ =\ \bar{a} \cos \tau_n  \bar{x} & =&\bar{a} \cos \tau_n  (0.5-x)
\eean
for $-0.25 \le \bar{x} \le 0.25$.
Since $\theta_n(x)$ is continuous at $x=\zeta_1=0.25$ which corresponds $\bar{x}=0.25$ we get the condition
\bean
\cos \tau_n  {\pi \over 4}&=& \bar{a} \cos \tau_n {\pi \over 4}
\eean
so we conclude that $ \bar{a}=1$.   

The  derivative jumps at $x=\zeta_1=0.25$  so we have
\bean
2\tau_n\sin {\tau_n  \over 4} &=&  {\sqrt{2}\over 2} \cos{\tau_n \over 4}
\eean
This leads to the equation
\bean
{\tau_n \over 4}&=& {\sqrt{2}\over 16} \cot{\tau_n \over 4}
\eean
Let $\sigma ={\tau_n \over 4}$ then we need to solve
\bean
\sigma_n &=&  {\sqrt{2}\over 16} \cot \sigma_n
\eean
There is exactly one root $ \sigma_n$ of this equation between $(n-1)\pi$ and
$(n-1/2)\pi$ for each $n=1,2,3,\ldots$ so 
\bean
4(n-1)<\tau_n< 4(n-1/2)
\eean
and so the eigenvalues $\eta_n$ of the error dynamics satisfy
\bean
-16(n-1)^2>\eta_n>-16(n-1/2)^2
\eean

We showed in Example 2 the closed loop eigenvalues under full state
feedback are $\mu_j=-(4j+2)^2\pi^2$ so this Linear Quadratic Gaussian
synthesis is asymptotically  stable.

\section{Conclusion}
We have explicitly derived an LQG synthesis for the heated/cooled rod under point actuation and point sensing.  The key is the completing 
the square technique.  We have already solved boundary control problems for the wave equation \cite{Kr21a}  and beam equation \cite{Kr21b} 
using the completing 
the square technique.  So we are confident that this technique can be used to solve the LQG Synthesis problems  for the wave and beam equations
under point actuation and point sensing.

\end{document}